\documentclass[Afour,sageh,times]{sagej}
\usepackage{tikz}
\usepackage{moreverb,url}
\usepackage[colorlinks,bookmarksopen,bookmarksnumbered,citecolor=red,urlcolor=red]{hyperref}

\newcommand\BibTeX{{\rmfamily B\kern-.05em \textsc{i\kern-.025em b}\kern-.08em
T\kern-.1667em\lower.7ex\hbox{E}\kern-.125emX}}

\def\volumeyear{2022}

\usepackage{multirow}
\usepackage{listings}
\def\Fb   {{\bf F}}
\def\ub   {{\bf u}}
\def\taub {\boldsymbol{\tau}}

\begin{document}

\runninghead{Loffeld et al.}

\title{Performance of explicit and IMEX MRI multirate methods on complex reactive flow problems within modern parallel adaptive structured grid frameworks}

\author{John J.~Loffeld\affilnum{1}, Andy Nonaka\affilnum{2}, Daniel R.~Reynolds\affilnum{3}, David J.~Gardner\affilnum{1}, Carol S.~Woodward\affilnum{1}}

\affiliation{\affilnum{1}Lawrence Livermore National Laboratory, USA\\
\affilnum{2}Lawrence Berkeley National Laboratory, USA\\
\affilnum{3}Southern Methodist University, USA
}

\corrauth{John J.~Loffeld, Lawrence Livermore National Laboratory, USA.}
\email{loffeld1@llnl.gov}

\begin{abstract}
Large-scale multiphysics simulations are computationally challenging due to the coupling of multiple processes with widely disparate time scales.  The advent of exascale computing systems exacerbates these challenges, since these systems enable ever increasing size and complexity.  In recent years, there has been renewed interest in developing multirate methods as a means to handle the large range of time scales, as these methods may afford greater accuracy and efficiency than more traditional approaches of using IMEX and low-order operator splitting schemes.  However, to date  there have been few performance studies that compare different classes of multirate integrators on complex application problems.  In this work, we study the performance of several newly developed multirate infinitesimal (MRI) methods, implemented in the SUNDIALS solver package, on two reacting flow model problems built on structured mesh frameworks.  The first model revisits the work of \cite{emmett2014high} on a compressible reacting flow problem with complex chemistry that is implemented using BoxLib but where we now include comparisons between a new explicit MRI scheme with the multirate spectral deferred correction (SDC) methods in the original paper.  The second problem uses the same complex chemistry as the first problem, combined with a simplified flow model, but run at a large spatial scale where explicit methods become infeasible due to stability constraints.  Two recently developed implicit-explicit (IMEX) MRI multirate methods are tested.  These methods rely on advanced features of the AMReX framework on which the model is built, such as multilevel grids and multilevel preconditioners.  The results from these two problems show that MRI multirate methods can offer significant performance benefits on complex multiphysics application problems and that these methods may be combined with advanced spatial discretization to compound the advantages of both.
\end{abstract}


\maketitle

\section{Introduction}
\label{sec:intro}

The advent of exascale computing will allow modeling of far larger and more complex systems than ever before.
In taking advantage of these systems, scientific simulations in areas such as climate, combustion, astrophysics, and fusion sciences are adding increasingly complex physical models with processes often evolving on widely varying time scales and requiring tighter, more accurate coupling.
For example, in atmospheric climate modeling, cloud microphysics, gas dynamics, and ocean currents each evolve at different rates.
In combustion, even for modest reaction mechanisms, the difference in time scales between stiff chemical kinetics, species/thermal diffusion, and acoustic wave propagation can easily span two or more orders of magnitude; convective time scales can be an additional one or two orders of magnitude slower.

Motivated by these large-scale multiphysics, multirate simulations composed of coupled models for interacting physical processes, we consider ordinary differential equation (ODE) initial value problems (IVPs) containing two time scales:
\begin{align}
\label{eq:IVP_fastslow2}
  y'(t) &= f^F(t,y) + f^S(t,y), \; t\in[t_0,t_f],\\
  y(t_0) &= y_0.\notag
\end{align}
Here $t$ is time, the dependent variable is $y \in \mathbb{R}^N$, $y'$ denotes $\mathrm dy/\mathrm dt$, $f^F: \mathbb{R} \times \mathbb{R}^{N} \to \mathbb{R}^N$ and $f^S: \mathbb{R} \times \mathbb{R}^{N} \to \mathbb{R}^N$ are partitions of the operator that advance at fast and slow timescales, respectively, and $y_0 \in \mathbb{R}^N$ is the initial state evolved from time $t_0$ to $t_f$.
For such problems, single-rate implicit time integration methods can become prohibitively expensive due to the cost and complexity of large, monolithic algebraic solvers.  Similarly, single-rate explicit integration methods may become infeasible due to stability constraints arising from the fastest process in the system.
Thus, practitioners frequently evolve the fast portion using small time steps, $h$, and the slow portion with a larger time step, $H \gg h$.
In many multiphysics circumstances, even these slow processes are frequently decomposed further into those for which explicit time-stepping methods provide sufficient accuracy and efficiency and those for which implicit methods are preferred to overcome restrictive stability limitations on the time step size for explicit approaches.  Such combinations are often treated using a mixed implicit-explicit (IMEX) scheme.  Thus, in this work we also consider the additional splitting of the problem into three components:
\begin{align}
  \label{eq:IVP_two_rate}
  y'(t) &= f^F(t,y) + f^I(t,y) + f^E(t,y), \; t\in[t_0,t_f],\\
  y(t_0) &= y_0.\notag
\end{align}
This may be achieved by splitting $f^S$ from \eqref{eq:IVP_fastslow2} into $f^I+f^E$, where the slow, nonstiff terms are grouped into $f^E$ and evolved using explicit methods, while the slow, stiff terms are grouped into $f^I$ and evolved using implicit methods.
Returning to the combustion example above, one could apply either of the multirate splittings \eqref{eq:IVP_fastslow2} or \eqref{eq:IVP_two_rate}, placing the stiff chemical kinetics in the operator $f^F$ and the remainder in $f^S$, or further partitioning the species/thermal diffusion into $f^I$ and acoustic wave propagation into $f^E$.

Traditional multirate approaches based on operator splitting have been used for many years for both \eqref{eq:IVP_fastslow2} and \eqref{eq:IVP_two_rate} with varying degrees of success, but these approaches have historically been limited to second-order.
In recent years, new multirate methods have been developed to allow for high-order integration algorithms that advance each partition at its own timescale potentially with a method tailored to each process.
These state-of-the-art multiphysics coupling schemes can further accelerate the discovery process by allowing even more efficient use of supercomputing resources.

In large scale simulation, the efficiency of the spatial discretization matters greatly as well.
A wide range of approaches for the spatial discretization of partial differential equations (PDEs) are based on structured grids.
Structured grid frameworks have a number of desirable properties in general and are increasingly favorable when leveraging exascale computing systems. Due to their regular grid structures, they offer a high degree of control over domain decomposition and parallel distribution of grids as well as low memory overhead due to simpler connectivity.
There is a rich history of numerical methods for structured grids that offer regular stencils and desirable data locality properties.  Regular stencils lead to efficient vectorization on manycore/GPU architectures, and improved data locality leads to fewer cache misses.
Block-structured adaptive mesh refinement (AMR) is a widely used strategy within the structured grid community that focuses computational power in regions of interest while retaining complete control over the grid structure. There are numerous applications and numerical methods that can take advantage of such a strategy in the fields of fluid mechanics, combustion, astrophysics, electromagnetism, and cosmology \citep{day2000numerical,fan2019maestroex,vay2018warp,almgren2013nyx}, just to name a few.

In addition to novel algorithms, newly developed numerical software targeted at efficient utilization of next generation supercomputing architectures is becoming increasingly available.
These packages are often open source, developed for portability, and can provide highly performant numerical methods within a structure that allows for continued development and extension of capabilities.
Both the United States Department of Energy and the European Union have made considerable investments in such software, and many scientific applications employ them 
on high performance systems.
In addition, the Extreme-scale Software Development Kit (xSDK) has developed community policies and standards that encourage such packages to adopt interoperability on many levels, both strengthening the capabilities available as well as their usability within complex scientific application codes \citep{xsdk:homepage}.

The goal of this paper is to understand the relative potential benefits of high-order multirate time integration schemes posed on adaptive, structured meshes.  Both the multirate methods and the structured AMR algorithms are instantiated within numerical software packages contained within the xSDK.
We consider reactive flow problems because they exhibit properties common to many multiphysics systems.  In particular, they include both fast (chemistry reaction systems) and slow (diffusion, acoustic, and/or convective) time dynamics and require a mixture of time integration methods for optimal performance.
A recent development in reactive flow simulation is the use of spectral deferred corrections (SDC) as a multirate temporal integrator. High-order multirate infinitesimal step (MRI) approaches developed over the last decade also offer potential speedups and robust alternatives to SDC on these problems.

There have been few performance studies for multirate integration on complex application problems.  Here we compare several MRI schemes supported by the ARKODE package \citep{reynolds2022arkode} in the SUNDIALS library \citep{hindmarsh2005sundials,gardner2020enabling} of time integrators and nonlinear solvers  to an existing SDC implementation for a combustion problem containing detailed kinetics and transport using the BoxLib block-structured AMR framework \citep{zhang2016boxlib}.
Additionally, we look at the performance of newly developed explicit and implicit-explicit (IMEX) multirate methods on a similar problem implemented with the AMReX block-structured AMR framework \citep{zhang2019amrex,zhang2020amrex}.
BoxLib is the predecessor to AMReX and contains a subset of features offered in AMReX (such as grid generation and MPI communication) and includes a complete SDC implementation of the combustion problem in our comparisons.
AMReX is a software framework developed by the Exascale Computing Project (ECP) Block Structured Adaptive Mesh Refinement Co-Design center and is designed to support structured adaptive grid calculations on massively parallel manycore/GPU machines. Also as part of the ECP, SUNDIALS has recently added several additional capabilities for integrating problems on GPU-based systems \citep{balos2021enabling}.
This work compares the performance of multirate methods for combustion using SDC and integrators available within the SUNDIALS time integrator library.
Our numerical experiments are performed using the AMReX and BoxLib frameworks.

The rest of the paper is organized as follows. The next section overviews the multirate and SDC time integration approaches and specific methods of interest to this work. Section \ref{sec:models} presents the mathematical formulation of the combustion problems we use to evaluate the performance of the integration methods. In section \ref{sec:impl} we provide overviews of AMReX and SUNDIALS and describe how we implemented our models in the combined SUNDIALS-AMReX framework.
In the next sections, we present results on two reacting flow problems.  Section \ref{sec:problem1} examines various explicit multirate schemes using a detailed combustion model, and section \ref{sec:problem2} presents results from an AMReX implementation of various IMEX multirate schemes for a modified reacting flow model. Finally, section \ref{sec:conclusions} summarizes the work and our findings.

\section{Multirate Time Integration}
\label{sec:multirate}

Many applications handle the varying time scales found in multiphysics systems using multirate approaches, where smaller time steps are used for rapidly evolving operators and larger steps are used for more slowly evolving components.
A popular, long-used approach to many multiphysics problems is use of low-order operator splitting with subcycling, e.g., first-order Lie--Trotter splitting \citep{mclachlanSplittingMethods2002} or second-order Strang--Marchuk splitting \citep{Marchuk1968,strangConstructionComparisonDifference1968a}.
However, these methods can suffer from a lack of stability and accuracy \citep{knoll2003balanced,estep2008posteriori}.
Higher order splitting methods have been developed, but they often require backward time integration \citep{goldman1996n,cervi2019comparison} which increases both the run time cost as well as the development time for code infrastructure.

Another highly flexible approach that multiphysics applications have used to handle multiple time scales is the use of SDC methods \citep{dutt2000spectral}.  In these methods, a series of correction equations are integrated to reduce error at each iteration in the series.  Typically, a low order method is used to integrate the system at each iteration, and these solutions are applied iteratively to build schemes of arbitrarily high accuracy.  The methods have been extended to handle solutions of multirate partial differential equation systems allowing for custom methods to be applied for each specific time scale and/or operator \citep{bourlioux2003high,layton2004conservative}.
While the SDC methods provide considerable freedom in the schemes used to integrate each scale and operator, they do require an integration over each point in time for each iteration, generally resulting in several passes over the time domain to generate a high order solution.  However, SDC methods have been used in the multiphysics community for some time \citep{huang2020high,nonaka2012deferred,pazner2016high,zingale2019improved} because of their relative ease in implementation by adapting established operator splitting implementations.
Specifically, high-order multirate approaches based on SDC have been shown to be effective on combustion problems with complex chemistry \citep{emmett2014high,emmett2019fourth}.
Within the combustion community, the long-time standard has been to employ high-order, explicit Runge-Kutta approaches such as those used in the S3D code \citep{chen2009terascale}.
In \cite{emmett2014high}, a multirate SDC approach for combustion was introduced in the ``SMC'' code, and a direct comparison between SMC and S3D demonstrated that integration in SMC can offer at least a 4x speedup over the widely used fourth-order explicit Runge-Kutta method for moderately-sized chemical mechanisms.

In the meantime, novel classes of multirate methods that show strong potential for a more efficient approach to high temporal order have been developed.
Starting with \cite{gear1984multirate} and continuing with
\cite{sand1992stability, gunther2016multirate-b, roberts2021implicit} and \cite{luan2020new}, multirate time integration methods were developed to address efficiency and accuracy issues associated with operator splitting approaches.
One particular multirate method, the multirate infinitesimal generalized additive Runge--Kutta (MRI-GARK) method \citep{sandu2019class}, has some distinct advantages over other multirate methods.  These schemes are a recent addition to the broader class of infinitesimal methods, including the multirate infinitesimal step (MIS) methods \citep{schlegel2011class, schlegel2009multirate}, that pioneered the approach.  Both MIS and MRI methods integrate a multirate problem of the form \eqref{eq:IVP_fastslow2} for one slow time step by alternating between evaluations of $f^S$ and advancing with $f^F$ for several fast steps.  With methods of this form, no part of the time step is repeated, and the method used for integrating the fast partition is left unspecified, allowing for specialized schemes tailored to the fast operator.  Thus MRI- and MIS-based methods allow extreme flexibility regarding methods applied to problem components, while simultaneously enabling high accuracy, without the requirement for backward-in-time integration or repeated evaluation of the step.

In the following subsections, we provide more detail about the specific classes of MRI and SDC approaches used in this work. 

\subsection{IMEX-MRI-GARK}
\label{sec:imex-mri}

The first class of multirate schemes that we consider is the implicit-explicit multirate infinitesimal GARK (IMEX-MRI-GARK) methods from \cite{chinomona2021implicitexplicit}.
These schemes are a recent addition to the broader class of infinitesimal methods, that generalize the previously-mentioned MIS and MRI-GARK methods constructed for problems of the form (\ref{eq:IVP_fastslow2}).
Methods within this family share two fundamental traits. First, they utilize a Runge--Kutta like structure for each time step $t_n\to t_n+H$, wherein the overall time step is achieved through computation of $s$ internal stages that are eventually combined to form the overall time step solution $y_{n+1}\approx y(t_n+H)$.
These internal stages roughly correspond to approximations $Y_i \approx y(t_n+c_i H)$, $i=1,\ldots,s$, where a specific scheme includes an array of coefficients $c\in\mathbb{R}^s$ that define these stage times. Second, these methods couple the fast and slow time scales by solving a set of \emph{modified} fast IVPs to obtain the internal stages.  Specifically, to compute the updated solution $y_{n+1}$ from $y_n$, these methods perform the algorithm:
\begin{enumerate}
\item Set $Y_1 = y_n$ and $T_1=t_n$.
\item For $i=2,\ldots,s$:
   \begin{itemize}
   \item[a.] Let $v(T_{i-1}) = Y_{i-1}$, $T_i = t_n+c_i H$ and $\Delta c_i=c_i-c_{i-1}$, and define a ``forcing'' polynomial $r_i(\theta)$ constructed from evaluations of $f^S(T_j,Y_j)$ where $j=1,\ldots,i$.
   \item[b.] Solve the ``fast'' IVP over $\theta\in[T_{i-1},T_i]$:
   \begin{align}
   \label{eq:modified_fast_IVP}
       v'(\theta) &= f^F(\theta, v(\theta)) + r_i(\theta).
   \end{align}
   \item[c.] Set $Y_i=v(T_i)$.
   \end{itemize}
\item Set $y_{n+1} = Y_{s}$.
\end{enumerate}
When $r_i(\theta)$ depends on $f^S(T_i,Y_i)$ the step 2b is implicit in $Y_i$; however when $\Delta c_i=0$, that step simplifies to a standard Runge--Kutta update.  Thus, infinitesimal methods that allow implicit treatment of $f^S$ typically only allow implicit dependence on $Y_i$ for stages where $\Delta c_i=0$, resulting in a ``solve-decoupled'' approach.
Here, the stages alternate between a standard diagonally-implicit slow nonlinear system solve and evolving the modified fast IVP. We note that since the step 2b may be solved using \emph{any} sufficiently-accurate algorithm, these methods offer extreme flexibility for practitioners to employ problem-specific integrators.

Within this infinitesimal family, the MIS methods first introduced by \cite{schlegel2009multirate} and \cite{schlegel2011class} treat the slow time scale operator $f^S$ explicitly, and result in up to $\mathcal{O}(H^3)$ accuracy.  The MRI-GARK methods introduced by \cite{sandu2019class} treat $f^S$ using either an explicit or a solve-decoupled implicit scheme, and result in up to $\mathcal{O}(H^4)$ accuracy.  The IMEX-MRI-GARK methods proposed by \cite{chinomona2021implicitexplicit} treat the slow time scale using ImEx methods, $f^S=f^I+f^E$, and result in up to $\mathcal{O}(H^4)$ accuracy.  Additionally, within IMEX-MRI-GARK methods, either of the components, $f^I$ or $f^E$, may be zero, resulting in an MRI-GARK method.

We note that these flexible infinitesimal formulations are not new, as the longstanding Lie-Trotter and Strang-Marchuk operator splitting methods similarly allow arbitrary algorithms at the fast time scale, as well as explicit, implicit, or IMEX treatment of the slow time scale.  However, both of these approaches use the forcing function $r_i(\theta)=0$, resulting in $\mathcal{O}(H)$ and $\mathcal{O}(H^2)$ accuracy, respectively.

\subsection{Multirate SDC Methods}
\label{sec:sdc}

Spectral deferred correction algorithms are a class of numerical
methods that represent the solution as an integral in time and iteratively solve a series
of correction equations designed to reduce the error.
The SDC approach is introduced in \cite{dutt2000spectral} for the solution of ordinary differential equations.
The correction equations are typically formed using a first-order time-integration scheme, but are applied iteratively to construct schemes of arbitrarily high accuracy.
In practice, the time step is divided into subintervals, and the solution at the end of each subinterval is iteratively updated using (typically) forward or backward Euler temporal discretizations.
Each set of correction sweeps over all the subintervals increases the overall order of the method by one, up to the underlying order of the numerical quadrature rule over the subintervals.
Thus, it is advantageous to define the subintervals using high-order numerical quadrature rules, such as Gauss-Lobatto, Gauss-Radau, Gauss-Legendre, or Clenshaw-Curtis.

There have been a number of works that have extended the SDC approach to handle physical systems with multiple disparate time scales.
\cite{minion2003semi} introduces a semi-implicit version for ODEs with stiff and non-stiff processes, such as advection–diffusion systems.
The correction equations for the non-stiff terms are discretized explicitly, whereas the stiff term corrections are treated implicitly.
\cite{bourlioux2003high} introduces a multirate SDC (MRSDC) approach for PDEs with advection-diffusion-reaction processes where advection terms are evaluated explicitly, reaction and diffusion terms treated implicitly, and different time steps are used for each process.
For example, the ``slow'' advection process can be discretized using Gauss-Lobatto quadrature over the entire time step while the stiffer diffusion/reaction processes can be sub-divided into Gauss-Lobatto quadrature(s) within each of the advection nodes.
This approach was also used successfully by \cite{layton2004conservative} for a conservative formulation of reacting compressible gas dynamics.
More recently, approaches have been developed for fourth-order, finite volume, adaptive mesh simulations \citep{emmett2019fourth}, as well as eighth-order finite differences \citep{emmett2014high}.
The latter paper, known as the ``SMC'' algorithm, forms the basis of the comparisons in this paper to the aforementioned MRI integrators.

Multirate SDC offers flexibility in (i) the number of physical processes and how they are grouped together, (ii) the choices of intermediate quadrature points each process uses, and (iii) the form of the correction equation for each process.
A long review can be found in the references above and mentioned within.
Here, we will describe one particular approach taken in the SMC algorithm that forms the basis of our comparisons that follow.

\begin{figure}[tb!]
\centering
\includegraphics[width=0.45\textwidth]{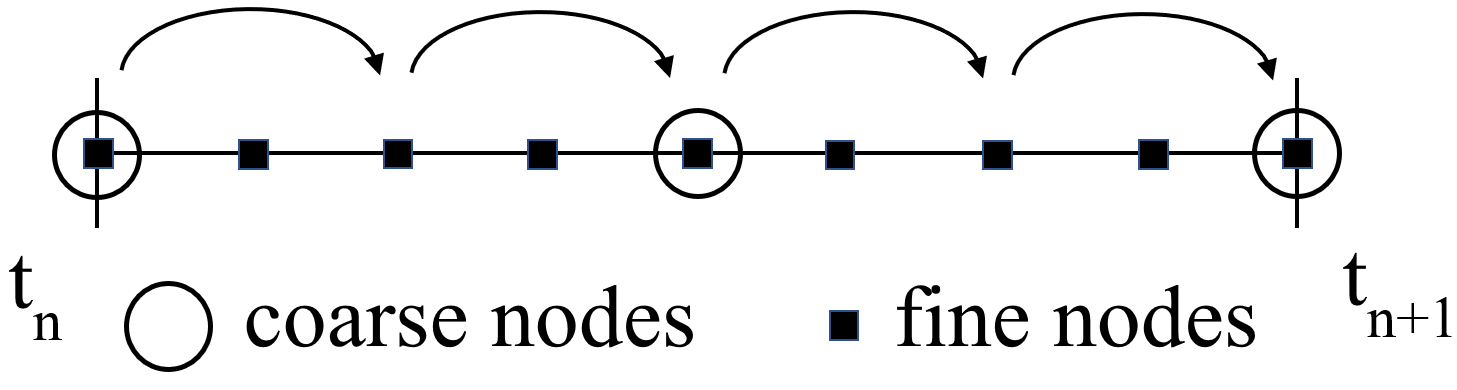}
\caption{Schematic of MRSDC-332 using the 3-point Gauss-Lobatto quadrature for the coarse nodes and applying a 3-point Gauss-Lobatto quadrature 2 times within each coarse subinterval. Note that the arrows demarcate the subintervals within the coarse interval.}
\label{fig:nodes}
\end{figure}

For the MRSDC methods in this paper we consider Gauss-Lobatto quadrature, noting that both endpoints of any given interval are included in the quadrature.
We divide the overall time step into a series of coarse substeps using nodes at time $t_p$, and each coarse substep is further divided into a series of finer substeps using nodes at time $t_q$.
With this choice of quadrature each coarse node at time $t_p$ corresponds to a location where a fine node also exists.
For comparison against the MRI schemes in this paper, we group advection and diffusion together as a ``slow'' explicit process, $f^S$, and evaluate these on the coarse nodes. We treat reactions as the ``fast'' explicit process, $f^F$, and evaluate their terms on fine nodes.
Other configurations are possible e.g., mixed implicit-explicit methods for the slow time scale with diffusion treated implicitly in $f^I$ \citep{bourlioux2003high}.
To describe the temporal scheme, we use the terminology ``MRSDC-XYZ'', where X is the number of coarse nodes, and we apply a Y-node quadrature Z times in each coarse subinterval.
As an example, Figure \ref{fig:nodes} shows a graphical illustration of an MRSDC-332 scheme with the 3-point Gauss-Lobatto quadrature, which consists of 3 equally-spaced points. Note that in general, there are $n_q = X + (X-1)(Z-1) + (X-1)(Y-2)Z$ total fine nodes.

The MRSDC scheme begins by supplying an initial guess for the solution at all nodes at time $t_q$
(typically the solution at $t_n$ is used), denoted $y_q^k$, where $k$ is the iteration index.
Each MRSDC iteration consists of a sweep over all the fine nodes at time $t_q$
in order, updating the solution from $y_q^k$ to $y_q^{k+1}$.
Each iteration in $k$ increases the overall order of accuracy by one (when using forward or backward Euler), up to the order of accuracy of the underlying quadrature.
When all processes are treated explicitly, the forward-Euler correction equation takes the form
\begin{eqnarray}
y_{q+1}^{k+1} = y_q^{k+1} &+& H_q\left(f^S(t_p,y_{p}^{k+1}) - f^S(t_p.y_{p}^k)\right)\nonumber\\
&+& H_q\left(f^F(t_q,y_{q}^{k+1}) - f^F(t_q.y_{q}^k)\right)\nonumber\\
&+& H_q S_{F}^{k,q} + H_q S_{S}^{k,q},\label{eq:SDC}
\end{eqnarray}
where $H_q = t_{q+1}-t_q$, and $p$
is the index of the coarse node closest to the left of time $t_q$,
$S_{F}^{k,q}$ is the numerical integral over $[t_q,t_{q+1}]$ of $f^F(t,y^k)$ using values from the fine nodes, and
$S_{C}^{k,q}$ is the numerical integral over $[t_q,t_{q+1}]$ of $f^S(t,y^k)$ using values from the coarse nodes.
In different formulations, some or all of the processes can be treated implicitly; for example, replacing $f^F(t_q,y_{q}^{k+1}) - f^F(t_q.y_{q}^k)$ with $f^F(t_{q+1},y_{q+1}^{k+1}) - f^F(t_{q+1}.y_{q+1}^k)$ in (\ref{eq:SDC}) represents a backward-Euler discretization of $f^F$,
forming an implicit update.
Equation (\ref{eq:SDC}) is applied over all $n_q$ nodes (for $q>0$)  a total of $N_{\rm sweep}$ times; for this explicit approach the computational expense is largely determined by the cost of the evaluations of $f^S$ and $f^F$.
The integration matrices $S_F$ and $S_S$ are (generally inexpensive) numerical integral quadratures over a small number of points.
For implicit or semi-implicit updates, the cost of each application of (\ref{eq:SDC}) can be much more expensive than the explicit case.

We thus see that both MRI and SDC approaches achieve high-order accuracy by coupling different (possibly multirate) physical processes via forcing terms.
Due to the iterative nature of SDC approaches, the computational expense tends to be larger than MRI methods of the same order, since MRI approaches advance the solution through a series of stages in ``one shot''.
We demonstrate the increased computational efficiency of MRI approaches in our results.
Also, it is worth noting that this difference in coupling strategies currently limits MRI approaches to at-best fourth order in time, whereas there are ostensibly no upper bounds to the temporal order of accuracy for SDC.

\section{Model Equations}
\label{sec:models}

We test the efficacy of the above multirate methods on two reacting flow model problems posed on structured meshes.
The first model is from the study by \cite{emmett2014high} which demonstrated that multirate SDC methods offer a significant performance improvement over single-rate SDC methods on reacting flow problems.
This model is similar to ones which have been successfully integrated using second-order splitting methods with reaction sub-cycling \citep{day2000numerical}, and the study by \cite{emmett2014high} shows that the benefit of multirate methods extends to higher order.  Re-using this same test problem allows the MRI methods to be tested under similar conditions and allows direct performance comparisons with the SDC methods.

The system of equations for this model uses the multicomponent reacting compressible Navier-Stokes equations given as
\begin{subequations}\label{eq:problem1}
\begin{eqnarray}
\frac{\partial\rho}{\partial t} &=& -\nabla\cdot(\rho\ub) ,\\
\frac{\partial\rho Y_k}{\partial t} &=& -\nabla\cdot(\rho Y_k\ub) + \rho\,\dot\omega_k - \nabla\cdot\Fb_k, \label{eq:smc_massfrac}\\
\frac{\partial\rho\ub}{\partial t} &=& -\nabla\cdot(\rho\ub\ub) - \nabla p + \nabla\cdot\taub, \\
\frac{\partial\rho E}{\partial t} &=& -\nabla\cdot[(\rho E+p)\ub] + \nabla\cdot\lambda\nabla T\nonumber\\
&&+ \nabla\cdot(\taub\cdot\ub) - \nabla\cdot\sum_k\Fb_k h_k.
\end{eqnarray}
\end{subequations}
Here, $\rho$ is the density, $\ub$ is the velocity, $p$ is the pressure, $E$ is the specific energy density (kinetic plus internal), $T$ is the temperature, $\taub$ is the viscous stress tensor, and $\lambda$ is the thermal conductivity.
For each of the chemical species $k$, $Y_k$ is the mass fraction, $\Fb_k$ is the diffusive flux, $h_k$ is the specific enthalpy, and $\dot\omega_k$ is the production rate.
The system is closed by an ideal gas mixture equation of state that specifies $p$ as a function of $\rho,T$, and $Y_k$,
\begin{equation}
p = \rho\mathcal{R}T\sum_j\frac{Y_j}{W_j},
\end{equation}
where $\mathcal{R}$ is the universal gas constant, and $W_j$ is the molecular weight of species $j$.
We use the standard full form of the viscous stress tensor,
\begin{equation}
\tau_{ij} = \eta\left(\frac{\partial u_i}{\partial x_j} + \frac{\partial u_j}{\partial x_i} - \frac{2}{3}\delta_{ij}\nabla\cdot\ub\right) + \xi\delta_{ij}\nabla\cdot\ub,
\end{equation}
where $\eta$ is the shear viscosity, and $\xi$ is the bulk viscosity.
The species diffusion flux is given by
\begin{equation}
\bar\Fb_k = \rho D_k\left(\nabla X_k + (X_k - Y_k)\frac{\nabla p}{p}\right),
\end{equation}
where $X_k$ is the mole fraction, and $D_k$ is the diffusion coefficient of species $k$.
In order to ensure mass is conserved, we impose a constraint that the sum of the species diffusion fluxes must equal zero.
We enforce this constraint using a correction velocity technique that modifies these fluxes with a mass-fraction weighted correction that ensures this constraint is satisfied, as commonly used in combustion algorithms and described in \cite{emmett2014high}.

In the second model, we compute the evolution of the mass fractions due to advection, diffusion, and reactions.
We use a simplified advective and diffusive transport model where we hold velocity, density, pressure, and the diffusion coefficients constant.
The governing equation for the mass fractions is
\begin{equation}
\frac{\partial Y_k}{\partial t} = \frac{\dot{\omega_k} W_k}{\rho} + \nabla\cdot (Y_k \ub) + \nabla\cdot D_k \nabla Y_k.\label{eq:problem2}
\end{equation}
We choose this model as it allows for testing of various MRI schemes using three-component mixed implicit-explicit systems
as in Equation (\ref{eq:IVP_two_rate}), and it serves as a proxy for the fully reacting flow problems we also consider.
The state-dependent diffusion coefficients in the first model require a more complex approach to implicit discretization that is beyond the scope of this paper.

\subsection{Spatial Discretization}
The spatial discretization for \eqref{eq:problem1} uses a finite difference approach with a uniform, single-level grid structure as described in \cite{emmett2014high}.  Centered, compact eighth-order spatial stencils are used in the interior of the domain, which require four cells on each side of a point where a derivative is computed.  As all calculations in our tests were done in parallel using spatial domain decomposition, four ghost cells are exchanged between each spatial partition when computing derivatives (see Section \ref{sec:block} below).  Stencils near non-periodic physical domain boundaries are handled as follows. Cells four away from a physical boundary use sixth-order centered stencils.  Cells three away use fourth-ordered centered stencils.  At cells two away, diffusion stencils use a centered second-order stencil while advection stencils use a biased third-order stencil.  For the boundary cell, diffusion uses a completely biased second-order stencil and advection a wholly biased third-order stencil.
While this mixture of discretizations theoretically results in overall second-order spatial accuracy, we replicate it here to allow for direct comparisons against the results from \cite{emmett2014high}.

The advection and diffusion operators of the second model, \eqref{eq:problem2}, are discretized using standard second-order finite-volume stencils.
These stencils require that one ghost cell be exchanged between spatial partitions during parallel computation of derivatives.  Test problems using this model all used periodic boundary conditions so the stencils remain second-order at physical boundaries.

\subsection{Temporal Discretization}
\label{sec:temporaldiscretization}
In this work, we explore a third-order explicit MRI method, as well as third- and fourth-order solve-decoupled IMEX-MRI-GARK methods.  The overall cost at the fast time scale for each method is commensurate with a single sweep over the time interval, while the overall cost at the slow time scale is commensurate with a standard additive Runge--Kutta IMEX (ARK-IMEX) method \citep{ascher1995implicit}.  We compare their performance on the two model problems against single- and multirate SDC methods, as well as a standard ARK-IMEX method.  The configuration and rough operation counts for a single step of each of the methods is as follows.
\begin{itemize}
    \item \textit{SDC}:  Explicit fourth-order single-rate SDC method that uses a 3-point quadrature and four iterations per step.
    \item \textit{MRSDC-338}:  Explicit fourth-order multirate SDC method reported in \cite{emmett2014high}.  Uses a 3-point quadrature for $f^S$, a 3-point quadrature repeated 8 times per iteration for $f^F$, and four iterations per step.
    \item \textit{MRSDC-352}:  Explicit fourth-order multirate SDC method reported in \cite{emmett2014high}.  Uses a 3-point quadrature for $f^S$, a high-order 5-point quadrature repeated only twice per iteration for $f^F$, and four iterations per step.
    \item \textit{MRI3}: The slow method is the third-order explicit method from \cite{knoth1998implicit} that requires four evaluations of $f^S$ and evolution of four modified fast IVPs \eqref{eq:modified_fast_IVP}.  At the fast time scale we use explicit Runge--Kutta methods of orders two through five, corresponding to Heun, \cite{bogacki19893} without use of the second-order embedding, \cite{zonneveld1963automatic}, and \cite{cash1990variable} respectively.
    \item \textit{ARK-IMEX}:  A fourth-order single-rate ARK-IMEX Runge-Kutta method from \cite{kennedy2003additive}.  Each step requires six evaluations of $f^E$ and the solution of five diagonally-implicit nonlinear systems using $f^I$.
    \item \textit{IMEX-MRI3}: The slow method is the IMEX-MRI-GARK3b method from \cite{chinomona2021implicitexplicit} and the fast method is Kutta’s third order method \cite{kutta1901}. Within each slow time step, this approach requires four evaluations of $f^E$, the solution of three diagonally-implicit nonlinear systems using $f^I$, and the evolution of three modified fast IVPs \eqref{eq:modified_fast_IVP}.
    \item \textit{IMEX-MRI4}: The slow method is the IMEX-MRI-GARK4 method from \cite{chinomona2021implicitexplicit} and the fast method is the classic fourth order Runge-Kutta method. Within each slow time step, this approach requires six evaluations of $f^E$, the solution of five diagonally-implicit nonlinear systems using $f^I$, and the evolution of five modified fast IVPs \eqref{eq:modified_fast_IVP}.
\end{itemize}


The model problem of \eqref{eq:problem1} is simulated using the multirate MRSDC-338, MRSDC-352, and MRI3 methods, plus the single-rate SDC.  The right hand side is kept unpartitioned for SDC but is divided into fast and slow partitions for the multirate methods (both MRSDC and MRI).  Since the simulation accuracy for all approaches hinges on resolution of the chemistry, the multirate methods put the chemical reaction term in Equation (\ref{eq:smc_massfrac}) into the fast partition,  $f^F = \rho \dot{\omega}_k$, and the remaining transport terms are put into the slow, $f^S$.  SDC must use the same step size for the transport as needed to accurately resolve the chemistry, whether that is necessary for them or not.

Problem \eqref{eq:problem2} is computed with the multirate IMEX-MRI3 and IMEX-MRI4 methods, as well as the standard ARK-IMEX method. For the same reason as the previous problem, the multirate integrators place the chemical reaction term into the fast partition, $f^F = \tfrac{\dot{\omega}_k W_k}{\rho}$.  The remaining slow terms are split into an explicitly stepped partition containing the advection term, $f^E = \nabla\cdot (Y_k \ub)$, and an implicitly stepped one with the diffusion term, $f^I = \nabla\cdot D_k \nabla Y_k$.  For ARK-IMEX, the chemistry and advection are assigned to the explicit partition, $f^E$, and the diffusion to the implicit one, $f^I$.  In that method, all partitions are stepped at the step size needed to resolve the chemistry, even if it is smaller than needed for the transport.

\section{Implementation}
\label{sec:impl}



Our implementation uses the AMReX \citep{zhang2019amrex,zhang2020amrex}
adaptive mesh refinement library and its predecessor package, BoxLib \citep{zhang2016boxlib} for spatial discretizations and the SUNDIALS library \citep{hindmarsh2005sundials,gardner2020enabling,reynolds2022arkode} for temporal discretizations.
Interfacing SUNDIALS with BoxLib allowed for an easier, direct comparison to previously developed code, and interfacing SUNDIALS with AMReX demonstrates the interoperability of the two libraries.

\subsection{Block structured meshes}
\label{sec:block}

\begin{figure}[tb!]
\centering
\includegraphics[width=0.25\textwidth]{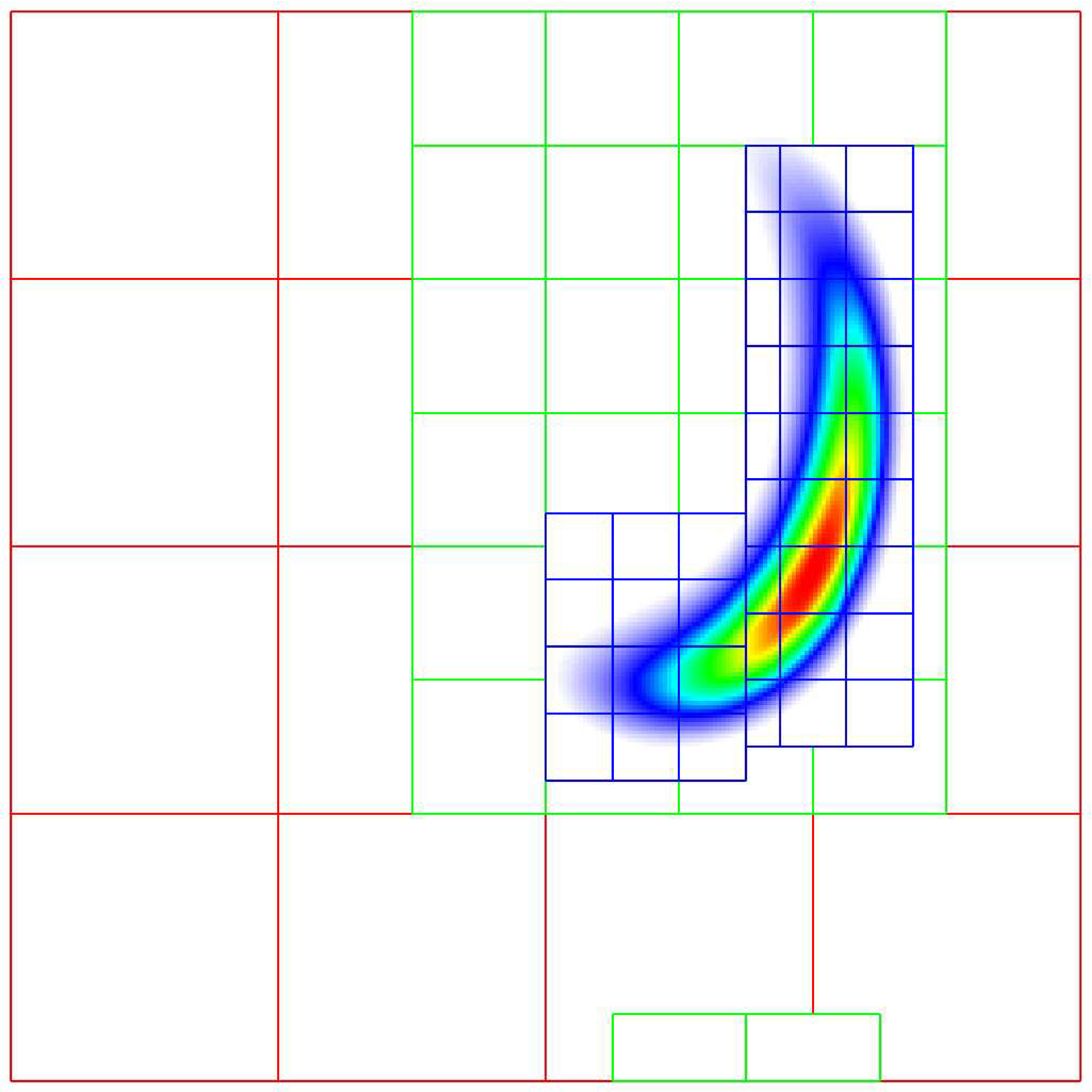}
\caption{In block-structured AMR, there is a hierarchy of logically rectangular grids. The computational domain on each AMR level is decomposed into a union of rectangular domains. In this figure there are 3 total levels.  The coarsest "level 0" grid covers the domain; bold red lines represent grid boundaries.
The green "level 1" grids contain cells that are a factor of two (or optionally four) finer than those at level 0.
The blue "level 2" grids contain cells that are a factor of two (or four) finer than those at level 1.
Note that there is no direct parent-child connection.}
\label{fig:grid}
\end{figure}

AMReX is a software framework that supports the development of block-structured AMR algorithms for solving systems of partial differential equations on current and emerging architectures.
AMReX supplies data containers and iterators for mesh-based fields, supporting data fields of any nodality.
There is no parent-child relationship between coarser grids and finer grids.
See Figure \ref{fig:grid} for a sample grid configuration.
In order to create refined grids, users define a ``tagging'' criteria to mark regions where refined grids are necessary.
The tagged cells are grouped into rectangular grids using the clustering algorithm in \cite{berger1991}.
There are user-defined parameters controlling the maximum length a grid can have, as well as forcing each side of the grid to be divisible by a particular factor.
Furthermore, grids can be subdivided if there is an excess of MPI processes with no load.
Proper nesting is enforced, such that a fine grid is strictly contained within the grids at the next-coarser level.

In this work we leverage (a) performance enhancements in the context of single-core and MPI performance, (b) built in grid generation, load balancing, multilevel communication patterns, (c) linear solvers via geometric multigrid, (d) efficient parallel I/O, and (d) visualization support through VisIt, ParaView, and yt.
The functionality in AMReX makes it easy to efficiently solve for the implicit terms in IMEX methods and is straightforward to use with both single-level and multi-level grids.
The code for the first model is implemented in the pure Fortran90 BoxLib framework.  For testing IMEX-MRI-GARK methods, the second model retains the complex chemistry from the first problem while taking advantage of the most updated version of AMReX.

AMReX uses an MPI+X strategy where grids are distributed over MPI ranks, and X represents a fine-grained parallelization strategy.
In recent years, the most common X in the MPI+X strategy has been OpenMP but AMReX now runs effectively on NVIDIA, AMD, and Intel GPUs.
The focus in this paper is on pure MPI; future efforts will be able to easily adopt the multicore/GPU strategies.
Load balancing over MPI ranks is done on a level-by-level basis using either options of a Morton-order space filling curve or knapsack algorithm.
AMReX efficiently manages and caches communication patterns between grids at the same level, grids at different levels, as well as data structures at the same level with different grid layouts.
The latter case is useful for when certain physical processes have widely varying computational requirements in different spatial regions, and grid generation with load balancing in mind may be weighted by this metric, rather than cell count.

\subsection{Time Integrators}
\label{sec:integrators}

The MRI methods evaluated in this work are provided by the SUNDIALS suite of time integrators and nonlinear solvers \citep{hindmarsh2005sundials, gardner2020enabling}. The suite is comprised of the ARKODE package of one-step methods for ODEs, the multistep integrators CVODE and IDA for ODEs and DAEs, respectively, the forward and adjoint sensitivity analysis enabled variants, CVODES and IDAS, and the nonlinear algebraic equation solver KINSOL. The ARKODE package provides an infrastructure for developing one-step methods for ODEs and includes time integration modules for stiff, non-stiff, mixed stiff/non-stiff, and multirate systems \citep{reynolds2022arkode}. ARKODE's MRIStep module implements explicit MIS methods and high-order explicit and solve-decoupled implicit and IMEX MRI-GARK methods. In this work, we solve each fast time scale ODE \eqref{eq:modified_fast_IVP} using ARKODE's ARKStep module which provides explicit, implicit, and ARK-IMEX Runge-Kutta methods, though a user may also provide their own fast scale integrator. The specific MRI methods evaluated in this work are built into ARKODE but users may alternatively supply their own coupling tables for MRI methods or Butcher tables for Runge-Kutta methods. At a minimum, to utilize the methods in ARKODE, an application must provide routines to evaluate the ODE right-hand side functions, and ARKODE must be able to perform vector operations on application data.

SUNDIALS packages utilize an object-oriented design and are written in terms of generic mathematical operations defined by abstract vector, matrix, and solver classes. Thus, the packages are agnostic to an application's particular choice of data structures/layout, parallelization, algebraic solvers, etc. as the specific details for carrying out the required operations are encapsulated within the corresponding classes. Several vector, matrix, and solver modules are provided with SUNDIALS, but users can also provide their own class implementations tailored to their use case. In this work, we supplied custom vector implementations to interface with the native BoxLib and AMReX data structures (one for BoxLib in Fortran, one each for the single and multi-level grid cases of AMReX in C++).  Often in-built functionality from BoxLib or AMReX could be used to aid these routines, but sometimes the cells in the grid needed to be iterated over and the operation computed per-cell explicitly.  The following is the listing for the vector multiplication operation (multiplies per-element like Matlab's \texttt{.*} operation) that leverages in-built AMReX operations in the multi-level case:
\begin{lstlisting}
void N_VProd_Multilevel(N_Vector x, N_Vector y, N_Vector z)
{
  int finest_level = NV_FINEST_LEVEL_M(x);

  // for each level in the grid
  for (int i = 0; i <= finest_level; i++)
  {
    // extract the level's AMReX MultiFabs
    // from the Sundials vectors
    MultiFab &mf_x = NV_MFAB_M(x,i);
    MultiFab &mf_y = NV_MFAB_M(y,i);
    MultiFab &mf_z = NV_MFAB_M(z,i);

    // get the number of components per cell
    sunindextype ncomp = mf_x.nComp();

    // do not include the ghost cells in the operation
    sunindextype nghost = 0;

    // operate on the AMReX MultiFabs
    MultiFab::Copy(mf_z, mf_x, 0, 0, ncomp, nghost);
    MultiFab::Multiply(mf_z, mf_y, 0, 0, ncomp, nghost);
  }
}
\end{lstlisting}

With explicit MRI methods, as used in the first problem \eqref{eq:problem1}, a vector implementing the necessary operations and routines for evaluating the right-hand side functions are all that is required for interfacing with SUNDIALS. For the IMEX-MRI methods, as used in the second problem \eqref{eq:problem2}, a nonlinear solver and a linear solver are needed to solve the nonlinear systems that arise at each implicit slow stage. Here we use the native SUNDIALS Newton solver and matrix-free GMRES implementation to solve these algebraic systems. As both solvers are written entirely in terms of vector operations no additional interfacing is required. However, to accelerate convergence of the GMRES iteration an additional wrapper function was created to leverage the AMReX multigrid linear solver as a preconditioner.
This solver is based on the ``MLABecLaplacian'' linear operator,
\begin{equation}
(\alpha I - \nabla\cdot\beta\nabla)\phi = \mathcal{R},
\end{equation}
where $\alpha$ and $\beta$ are spatially-varying coefficients.

\section{Explicit time stepping on compressible reacting flow}
\label{sec:problem1}

The results in this section are motivated by earlier success of the explicit MRSDC methods on the SMC combustion problem \citep{emmett2014high} in Equation \eqref{eq:problem1}.  Here we compare the performance of various MRI and MRSDC methods on this same model and discuss how their structural differences affect their relative performance.

\subsection{Problem setup}
\label{sec:problem1-setup}

The experiment in this section is close in configuration to the three-dimensional flame problem from Section 5.2.2 in \cite{emmett2014high}.  The problem domain is a cubic box with 1.6~cm on each side.  The boundaries are periodic in the y and z directions and inlet/outlet in the x-direction.
The chemistry mechanism for the flame is the GRI-MECH 3.0 network \citep{frenklach1995gri} that contains 53 species and 325 reactions. Combined with the other state variables such as density and velocity, the resulting system has 58 unknowns at each spatial location.
The initial conditions are taken from the PREMIX code from the CHEMKIN library \citep{kee1998premix}, where we use a precomputed realistic distribution of species in a premixed methane flame in 1D along the x-axis
with offsets in the pattern varying in the x-y plane by perturbing spherical harmonics, resulting in a bumpy pattern at the front.  The problem is integrated over a total time of $t=2.0\times 10^{-6}$ s using coarse time step sizes that range from $H=2.0\times 10^{-6}$ to $H=1.5625\times 10^{-8}$ s.

The domain interior is discretized with $32^3$ points.  Since there are $58$ unknowns per spatial location, this discretization gives a total of 1,900,544 unknowns over the interior of the domain.  The problem is parallelized by splitting each dimension in half, resulting in 8 boxes distributed over 8 MPI ranks.  Since the discretization also requires an additional four layers of ghost cells at each face, the total number of data in the problem becomes 6,414,336 values when including ghost cells.

As with \cite{emmett2014high}, we test this problem on the single-rate SDC and the two multirate SDC methods MRSDC-338 and MRSDC-352.  Here we also compare against MRI3.  As noted in Section \ref{sec:temporaldiscretization}, the MRI3 scheme is paired with fast methods of orders two through five.  For each method we perform k fast substeps per slow step (where that ratio is held constant over the problem), and denote the resulting combinations as MRI-P(k), where is P the order of the fast method.  Through exhaustive testing, on this application we found the following options to be the most efficient for each fast method order: MRI-2(8), MRI-3(4), MRI-4(8), and MRI-5(4).

The results were gathered on a node of the Quartz cluster at Lawrence Livermore National Laboratory.  Each node is composed of two Intel Xeon E5-2695 v4 chips.  The tests used 8 ranks of MPI parallelism distributed over the cores of the node.

\subsection{Results}
\label{sec:problem1-results}



The range of magnitudes of the different solution components (e.g. density, velocity, mass fractions of different species) is quite large in this problem, so a norm of error over the solution would not properly capture the total error over all components without using an appropriate weighting vector in the norm.  As in \cite{emmett2014high}, we find the order of convergence of the methods to not depend on which solution component is chosen.  Therefore, rather than use a norm over all components, we measure the error as the max norm in density.  The results are similar for other components.

The accuracy of the methods is tested against a reference solution calculated using a six-stage fourth-order explicit Runge-Kutta method \citep{kennedy2000low} with a step size of $5.0 \times 10^{-10}$ s.
The integration of each method is done using constant time stepping, starting each with the coarsest time step that exhibited numerical stability and halving the step for each successive run.
The coarsest allowable time step for each method can be seen in Table \ref{table:convergence_rates}.
Not included in the table is the single-rate SDC case, which required a step size no greater than $6.25 \times 10^{-8}$ s to be stable.
We note that the single-rate restriction of $\mathcal O(10^{-8})$ is consistent with the analysis of the eigenvalues in the Jacobian of this reaction network \citep{wartha2020characteristic}.

\begin{figure}[tb!]
\centering
\includegraphics[width=0.45\textwidth]{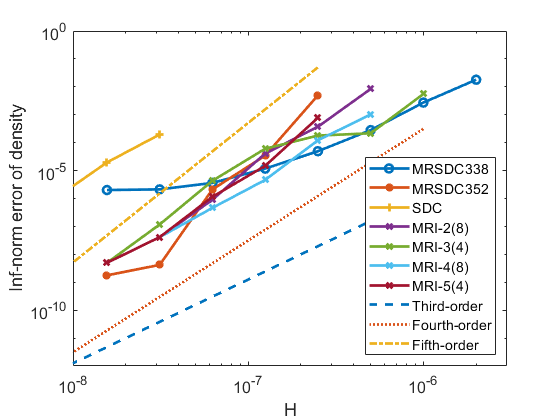}
\caption{Error in each method vs.~time step for the explicit methods from Section \ref{sec:problem1}.  
Each of the MRI methods exhibit at least third-order accuracy as expected. The convergence rates for the MRSDC approaches is highly non-monotonic.}
\label{fig:convergence}
\end{figure}

%
The rate of convergence for each of the methods is shown in Figure \ref{fig:convergence} and the same data is shown in tabular form in Table \ref{table:convergence_rates}.
We note that the convergence rates for all cases remains very similar if we consider the $L_2$ norm instead of $L_\infty$ norm.
In general, the MRSDC methods exhibit highly non-monotonic convergence rates, despite their theoretical fourth-order accuracy.
This problem could possibly be alleviated by further reducing the step size, or increasing the number of SDC iterations, as has been demonstrated in other works using SDC for reacting flow \citep{nonaka2012deferred,nonaka2018conservative,pazner2016high}.
The rate of convergence of the MRI methods is a balance between the accuracy of the fast and slow partitions.  The MRI-4 and MRI-5 methods are asymptotically limited to third-order accuracy due to the third-order accuracy of the slow integrator, however the slope of the convergence plot is better than four in some ranges when the reactions dominate the overall accuracy.  The MRI-2 method is asymptotically limited to only second-order as step sizes become small, but the slope of the line is at least three for the ranges of step sizes used in the tests.  Overall, the MRI methods exhibit convergence rates consistent with the expected order of accuracy, and appear to reach the asymptotic range of convergence at larger step sizes than the MRSDC approaches.

\setlength{\tabcolsep}{5pt} 
\begin{table*}[tb!]
\caption{$L_\infty$ error and convergence rates for density using various multirate approaches on a methane flame test from Section \ref{sec:problem1}.  The omitted values for largest time steps indicate simulations that became unstable and were unable to finish.
}
\label{table:convergence_rates}
\footnotesize{
\begin{tabular}{lcccccccccccc}
Method & \multicolumn{2}{c}{MRSDC338} & \multicolumn{2}{c}{MRSDC352} & \multicolumn{2}{c}{MRI-2(8)} & \multicolumn{2}{c}{MRI-3(4)} & \multicolumn{2}{c}{MRI-4(8)} & \multicolumn{2}{c}{MRI-5(4)} \\
$H$ [s]&  $L_\infty$ error &  rate &  $L_\infty$ error &  rate &  $L_\infty$ error &  rate &  $L_\infty$ error &  rate &  $L_\infty$ error &  rate &  $L_\infty$ error &  rate \\
\hline
2.0$\times 10^{-6}$    & 1.80$\times 10^{-2}$ &      &                      &      &                      &      &                      &      &                      &      &                       & \\
1.0$\times 10^{-6}$    & 2.73$\times 10^{-3}$ & 2.72 &                      &      &                      &      & 5.50$\times 10^{-3}$ &      &                      &      &                       & \\
5.0$\times 10^{-7}$    & 2.81$\times 10^{-4}$ & 3.28 &                      &      & 8.40$\times 10^{-3}$ &      & 2.17$\times 10^{-4}$ & 4.48 & 1.03$\times 10^{-3}$ &      &                       & \\
2.5$\times 10^{-7}$    & 4.94$\times 10^{-5}$ & 2.51 & 4.95$\times 10^{-3}$ &      & 3.72$\times 10^{-4}$ & 4.50 & 1.79$\times 10^{-4}$ & 0.27 & 1.21$\times 10^{-4}$ & 3.09 & 7.77$\times 10^{-4}$ & \\
1.25$\times 10^{-7}$   & 1.17$\times 10^{-5}$ & 2.07 & 3.44$\times 10^{-5}$ & 7.17 & 4.07$\times 10^{-5}$ & 3.20 & 6.07$\times 10^{-5}$ & 1.56 & 4.63$\times 10^{-6}$ & 4.71 & 1.48$\times 10^{-5}$ & 5.71 \\
6.25$\times 10^{-8}$   & 3.62$\times 10^{-6}$ & 1.69 & 2.12$\times 10^{-6}$ & 4.02 & 9.06$\times 10^{-7}$ & 5.89 & 4.31$\times 10^{-6}$ & 3.81 & 4.61$\times 10^{-7}$ & 3.33 & 1.15$\times 10^{-6}$ & 3.68 \\
3.125$\times 10^{-8}$  & 2.12$\times 10^{-6}$ & 0.77 & 4.19$\times 10^{-9}$ & 8.98 & 4.12$\times 10^{-8}$ & 4.46 & 1.18$\times 10^{-7}$ & 5.20 & 4.12$\times 10^{-8}$ & 3.48 & 4.12$\times 10^{-8}$ & 4.80 \\
1.5625$\times 10^{-8}$ & 2.01$\times 10^{-6}$ & 0.08 & 1.77$\times 10^{-9}$ & 1.24 & 5.02$\times 10^{-9}$ & 3.04 & 4.95$\times 10^{-9}$ & 4.57 & 4.96$\times 10^{-9}$ & 3.05 & 4.96$\times 10^{-9}$ & 3.05
\end{tabular}
}
\end{table*}
\setlength{\tabcolsep}{6pt} 

\begin{figure}[tb!]
\centering
\includegraphics[width=0.45\textwidth]{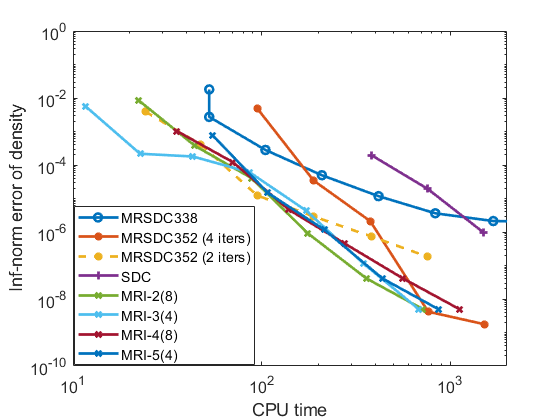}
\caption{Error in each method vs.~CPU time for the explicit methods from Section \ref{sec:problem1}.
}
\label{fig:basicperformance}
\end{figure}

Figure \ref{fig:basicperformance} shows an efficiency diagram comparing the performance of the various integrators, plotting the error in density as a function of wallclock time.  The figure confirms the results of \cite{emmett2014high} that the MRSDC methods have a significant performance advantage over the single-rate SDC method on this problem.
We also see that in general, the MRI methods produce smaller error for a given wallclock time compared to the SDC approaches.  The exception is for the smallest step sizes reported, the MRSDC352 4-iteration performs comparably to the MRI approaches.
As discussed in Section \ref{sec:sdc}, the performance benefits of MRI approaches arise because the MRI methods need fewer calls to the right-hand-side functions per step than the SDC methods.
To further quantify this point, Table \ref{table:scheme_compare_6.25e-8} presents the $H = 6.25 \times 10^{-8}$~s data from Figure \ref{fig:basicperformance} / Table \ref{table:convergence_rates}, including the number of function calls, runtime, and error using four different approaches.
The error for the MRI approaches is comparable to or smaller than the SDC approaches, in addition to offering runtimes that are notably smaller due to the reduced number of slow and fast evaluations.

\begin{table}[tb!]
\small\sf\centering
\caption{Scheme comparison using $H = 6.25\times 10^{-8}$ [s] over 32 time steps for the explicit methods from Section \ref{sec:problem1}.
}
\label{table:scheme_compare_6.25e-8}
\begin{tabular}{ccccc}
\toprule
Scheme    & \multicolumn{2}{c}{\# Evals} & Runtime   & Error                 \\
          & Slow         & Fast          & {[}$s${]} & {[}$g/cm^3${]}        \\
\midrule
MRSDC-338 & 129          & 7,169         & 833       & 3.62 $\times10^{-6}$  \\
MRSDC-352 & 257          & 2,049         & 376       & 2.12 $\times 10^{-6}$ \\
MRI-2(8)  & 97           & 993           & 176       & 9.06 $\times 10^{-7}$ \\
MRI-3(4)  & 97           & 929           & 173       & 4.31 $\times 10^{-6}$ \\
MRI-4(8)  & 97           & 1857          & 274       & 4.61 $\times 10^{-7}$ \\
MRI-5(4)  & 97           & 1345          & 216       & 1.15 $\times 10^{-6}$ \\
\bottomrule
\end{tabular}
\end{table}

\section{IMEX time stepping on a reacting flow problem implemented in AMReX}
\label{sec:problem2}

For IMEX multirate tests, we use the simplified model defined by \eqref{eq:problem2}.  The chemistry reaction mechanism is the same as for the compressible reacting flow problem of Section \ref{sec:problem1}.  However, here the diffusion term is much stiffer to represent the type of diffusive stiffness that must be confronted on large scale transport problems.  The advection term is non-stiff.  As with the previous problem, the chemistry uses a small timestep to resolve the fast dynamics and the transport uses coarse time steps, but now the diffusion term must be handled implicitly to deal with its stiffness.

\subsection{Problem setup}
\label{sec:problem2-setup}

The problem domain is a 3D cube $[-0.05, 0.05]$~cm on a side with periodic boundaries in each direction.  The chemistry mechanism is the same GRI-MECH 3.0 network \citep{frenklach1995gri} as the previous problem (53 species, 325 reactions).  The model is integrated over a total time of $t = 2.0 \times 10^{-6}$ s.

The problem is initialized using the same predetermined realistic distribution of species in 1D for a fuel mix as that used in the previous example.  However, here the 1D distribution is swept angularly around the origin with high concentration offset from the middle, resulting in a ring of fuel centered around the center.  The temperature is set as $1500$ K everywhere, the pressure is set at $1$ atmosphere everywhere, and the density is calculated from the standard equation of state from the pressure and temperature at each cell.  The velocity of the system is $\langle 0.001, 0.001, 0.001 \rangle$ cm/s uniformly across the domain.

The problem is tested using ARK-IMEX, IMEX-MRI3 and IMEX-MRI4.  Details of the methods and the partitioning of the model terms are described in Section \ref{sec:temporaldiscretization}.  The time stepping is done with constant step sizes, with ten fast steps taken for each slow step in the multirate methods.  The multirate methods are thus labeled as ``IMEX-MRI3-10'' and ``IMEX-MRI4-10'' in what follows.


\begin{table}[tb!]
\centering
\begin{tabular}{c|c}
\toprule
\multirow{ 2}{*}{$H$ [s]} & ARK-IMEX \\
                             &   tol   \\
\midrule
$1.0 \times 10^{-7}$   & $10^{-5}$     \\
$5.0 \times 10^{-8}$   & $10^{-7}$     \\
$2.5 \times 10^{-9}$   & $10^{-7}$     \\
$1.25 \times 10^{-9}$  & $10^{-8}$     \\
$6.25 \times 10^{-10}$ & $10^{-12}$    \\
\bottomrule
\end{tabular}
\caption{Solver tolerance values used at each step size for the ARK-IMEX method from Section \ref{sec:problem2}.}
\label{tab:tolerances_ao}
\end{table}

\begin{table}[tb!]
\centering
\begin{tabular}{c|c|c}
\toprule
\multirow{ 2}{*}{$H$ [s]}  & IMEX-MRI3-10 &  IMEX-MRI4-10 \\
                           & tol          &  tol     \\
\midrule
$2.0 \times 10^{-6}$   & $10^{-4}$                                              & $10^{-4}$                                             \\
$1.0 \times 10^{-6}$   & $10^{-5}$                                              & $10^{-6}$                                             \\
$5.0 \times 10^{-7}$   & $10^{-6}$                                              & $10^{-7}$                                             \\
$2.5 \times 10^{-7}$   & $10^{-8}$                                              & $10^{-7}$                                             \\
$1.25 \times 10^{-7}$  & $10^{-8}$                                              & $10^{-9}$                                             \\
$6.25 \times 10^{-8}$  & $10^{-12}$                                             & $10^{-11}$                                            \\
$3.125 \times 10^{-8}$ & $10^{-10}$                                             & $10^{-14}$                                            \\
\bottomrule
\end{tabular}
\caption{Solver tolerance values used at each step size for the IMEX-MRI3-10 and IMEX-MRI4-10 methods from Section \ref{sec:problem2}.}
\label{tab:tolerances_mri}
\end{table}

For the implicit partitions, the implicit stage solves are done using SUNDIALS' Newton iteration and  matrix-free GMRES solver.  No preconditioning is used.  To prevent iterations being terminated before reaching tolerance, the maximum number of nonlinear and linear iterations is relaxed from the default to 10 and 100 respectively.  However, solver tolerances also need to be set for the nonlinear and linear solvers in order to avoid ``over-solving'' the algebraic systems of equations as well.  The nonlinear solver tolerances were chosen on a per-method and per-timestep size basis.  The chosen solver tolerance values are listed in Tables \ref{tab:tolerances_ao} and \ref{tab:tolerances_mri}.  The reported values were chosen by sweeping each case over a range of tolerances from $10^{-1}$ to $10^{-15}$, with jumps of a factor of $10$ between each, with the reported value being the coarsest one for which both the H2 and H mass fraction error is perturbed by less than a factor of $0.01$.  The tolerance value is used as an absolute tolerance in the WRMS norm for the nonlinear solver and the 2-norm for the linear solver over the entire solution vector, i.e. over all species at once.  Safety factors of 0.1 and 0.05 are applied for the nonlinear and linear solvers, respectively.

The reference solution is computed using the fifth-order explicit Runge--Kutta method from \cite{cash1990variable} with a step size of $10^{-10}$, which is an order of magnitude smaller than the smallest fast step size used in all the runs.

The grid resolution is 128 cells on a side.  Diffusion coefficients are computed in the same manner and library code as the SMC problem, using temperatures, molar fractions, and specific heat at constant pressure.  The calculation of the diffusion coefficient is computationally expensive, and dominates the cost of the implicit terms more than the nonlinear solves for the tested configuration of the problem.  The coefficients are not updated within the solves, so the computational cost of that portion is not affected by the changes in iteration counts with stiffness.  However, in order to make the problem sufficiently stiff to warrant an implicit method, the computed coefficient in each cell is scaled by a factor of 1000.  The problem is thus not physically realistic, but is meant to exercise the same code components that are expensive in a reactive-flow problem at a level of stiffness that would be found in highly resolved realistic cases.  The scale of the diffusion coefficient combined with the grid resolution results in a problem that cannot be completed successfully using explicitly evolved diffusion at the coarser step sizes of the runs.

The problem is run with 512 total ranks over 16 nodes on the Quartz cluster at Lawrence Livermore National Laboratory.

\subsection{Results}
\label{sec:problem2-results}

\begin{figure}[tb!]
    \centering
	\includegraphics[width=0.45\textwidth]{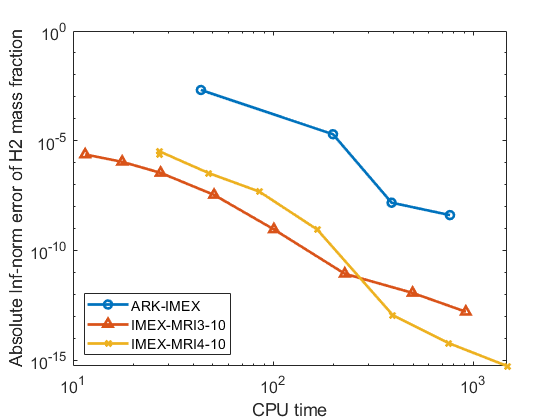}
	\caption{Efficiency diagrams for the IMEX methods from Section \ref{sec:problem2} with error measured for the H2 mass fraction.}
	\label{fig:amrex_efficiency_h2}
\end{figure}

The results of the runs are shown in Figure \ref{fig:amrex_efficiency_h2}. The error is shown for the species H2, but as with the SMC problem, the outcome is insensitive to which species is used to measure error.

We see in the plots that there is a significant performance advantage when using multirate over single-rate IMEX methods for this problem.  We tried the problem over a range of fast-to-slow step size ratios and found ten-to-one to be the best balance.  Due to the expensive diffusion coefficient calculation, reducing the number of coarse steps that must be computed per fast step results in a large time savings.

An interesting phenomenon from Tables \ref{tab:tolerances_ao} and \ref{tab:tolerances_mri} is that the solver tolerances needed by the multirate methods were more stringent than for the single-rate method for the same coarse step size.  This seems to be due to the increased accuracy of the MRI solutions at those step sizes, thereby necessitating tighter algebraic solver tolerances.

Table \ref{tab:scaling_tol} shows the scaling in CPU time and change in nonlinear and linear iterations for the problem as the grid resolution is refined.  The results are from using the IMEX-MRI3-10 method over eight time steps, i.e. at about mid-height in Figure \ref{fig:amrex_efficiency_h2}.  Even at a somewhat finer step size than the coarsest case, the problem is moderately stiff so the number of linear iterations grows with the grid resolution.  However, the primary cost is the computation of the diffusion coefficients, which are not done within the linear iterations.  As a result, the total cost goes up only slowly with stiffness within the range of resolutions tested.  If the problem were to continue to grow in size, eventually the solve cost would take over the cost of the diffusion coefficient calculation and a preconditioner would be needed.  As mentioned in Section \ref{sec:integrators}, AMReX includes the MLABecLaplacian multigrid solver that works over multiple AMR grid levels.  Here, we use this preconditioner for this single-level problem, with results comparing the un-preconditioned and preconditioned solvers in Tables \ref{tab:scaling_tol} and \ref{tab:precond_multi_scaling_tol}, respectively.  Due to the overhead in invoking the preconditioner, for coarse problems the un-preconditioned solve is more efficient, but once the grid reaches at least 128 cells per side the preconditioner becomes cost effective.  Tables \ref{tab:stiff_scaling_tol} and \ref{tab:precond_multi_stiff_scaling_tol} push this experiment further, by artificially increasing the problem stiffness through scaling the diffusion coefficients by another factor of 10.  Here again, preconditioning is effective at grid sizes exceeding 64 cells per side, and becomes critical as the grids are increased to 256 cells per side.

\begin{table}[tb!]
\centering
\begin{tabular}{l|l|l|l|ll}
\toprule
\multicolumn{1}{c|}{Cells}    & \multirow{2}{*}{Ranks} & \multicolumn{1}{|c}{Solver} & \multicolumn{1}{|c}{CPU}  & \multicolumn{1}{|c}{Nonlinear} & \multicolumn{1}{|c}{Linear}  \\
\multicolumn{1}{c|}{per side} &                        & \multicolumn{1}{|c}{tol}    & \multicolumn{1}{|c}{Time} & \multicolumn{1}{|c}{Iters.}    & \multicolumn{1}{|c}{Iters.} \\
\midrule
32   & 8     & $10^{-8}$   & 35.6   & \multicolumn{1}{l|}{48}  & 119    \\
64   & 64    & $10^{-8}$   & 46.5   & \multicolumn{1}{l|}{48}  & 142   \\
128  & 512   & $10^{-8}$   & 69.0   & \multicolumn{1}{l|}{48}  & 209   \\
256  & 4096  & $10^{-12}$  & 129.1  & \multicolumn{1}{l|}{72}  & 1150  \\
\bottomrule
\end{tabular}
\caption{The CPU time and the number of nonlinear and linear iterations required by the IMEX-MRI3-10 method when computing the problem from Section \ref{sec:problem2} over 8 steps as the grid resolution is refined and the stiffness increases.  The solver tolerance is chosen fine enough to not limit the accuracy of the calculations.}
\label{tab:scaling_tol}
\end{table}

\begin{table}[tb!]
\centering
\begin{tabular}{l|l|l|l|ll}
\toprule
\multicolumn{1}{c|}{Cells}    & \multirow{2}{*}{Ranks} & \multicolumn{1}{|c}{Solver} & \multicolumn{1}{|c}{CPU}  & \multicolumn{1}{|c}{Nonlinear} & \multicolumn{1}{|c}{Linear}  \\
\multicolumn{1}{c|}{per side} &                        & \multicolumn{1}{|c}{tol}    & \multicolumn{1}{|c}{Time} & \multicolumn{1}{|c}{Iters.}    & \multicolumn{1}{|c}{Iters.} \\
\midrule
32   & 8     & $10^{-8}$   & 45.92  & \multicolumn{1}{l|}{48}  & 48  \\
64   & 64    & $10^{-8}$   & 64.45  & \multicolumn{1}{l|}{48}  & 48  \\
128  & 512   & $10^{-8}$   & 66.63  & \multicolumn{1}{l|}{48}  & 48  \\
256  & 4096  & $10^{-12}$  & 112.1  & \multicolumn{1}{l|}{72}  & 60  \\
\bottomrule
\end{tabular}
\caption{The same experiment as in Table \ref{tab:scaling_tol}, but when using the AMReX MLABecLaplacian solver for the implicit Diffusion component.  For the grids with 32 and 64 cells per side, the preconditioner is a net disadvantage.  However, for the grid with 128 cells per side, which is the same as the main experiments for this section, the preconditioner starts to show a slight advantage.  For the grid that is 256 cells per side, the advantage widens.}
\label{tab:precond_multi_scaling_tol}
\end{table}

\begin{table}[tb!]
\centering
\begin{tabular}{l|l|l|l|ll}
\toprule
\multicolumn{1}{c|}{Cells}    & \multirow{2}{*}{Ranks} & \multicolumn{1}{|c}{Solver} & \multicolumn{1}{|c}{CPU}  & \multicolumn{1}{|c}{Nonlinear} & \multicolumn{1}{|c}{Linear}  \\
\multicolumn{1}{c|}{per side} &                        & \multicolumn{1}{|c}{tol}    & \multicolumn{1}{|c}{Time} & \multicolumn{1}{|c}{Iters.}    & \multicolumn{1}{|c}{Iters.} \\
\midrule
32   & 8     & $10^{-8}$   & 40.9   & \multicolumn{1}{l|}{48}  & 241  \\
64   & 64    & $10^{-8}$   & 64.6   & \multicolumn{1}{l|}{48}  & 439  \\
128  & 512   & $10^{-8}$   & 102.8  & \multicolumn{1}{l|}{48}  & 875  \\
256  & 4096  & $10^{-12}$  & 853.6  & \multicolumn{1}{l|}{72}  & 4796 \\
\bottomrule
\end{tabular}
\caption{The same experiment as in Table \ref{tab:scaling_tol}, without a preconditioner, but with a diffusion coefficient multipler of 1e4 instead of 1e3.}
\label{tab:stiff_scaling_tol}
\end{table}

\begin{table}[tb!]
\centering
\begin{tabular}{l|l|l|l|ll}
\toprule
\multicolumn{1}{c|}{Cells}    & \multirow{2}{*}{Ranks} & \multicolumn{1}{|c}{Solver} & \multicolumn{1}{|c}{CPU}  & \multicolumn{1}{|c}{Nonlinear} & \multicolumn{1}{|c}{Linear}  \\
\multicolumn{1}{c|}{per side} &                        & \multicolumn{1}{|c}{tol}    & \multicolumn{1}{|c}{Time} & \multicolumn{1}{|c}{Iters.}    & \multicolumn{1}{|c}{Iters.} \\
\midrule
32   & 8     & $10^{-8}$   & 48.0  & \multicolumn{1}{l|}{48}  & 48 \\
64   & 64    & $10^{-8}$   & 68.8  & \multicolumn{1}{l|}{48}  & 48 \\
128  & 512   & $10^{-8}$   & 82.3  & \multicolumn{1}{l|}{48}  & 48 \\
256  & 4096  & $10^{-12}$  & 145.8 & \multicolumn{1}{l|}{72}  & 72 \\
\bottomrule
\end{tabular}
\caption{The same experiment as in Table \ref{tab:stiff_scaling_tol}, but using the AMReX MLABecLaplacian solver for the implicit Diffusion component.  For two coarsest grids, using the preconditioner is still a net disadvantage, but for the two finest grids the preconditioner provides a significant advantage.}
\label{tab:precond_multi_stiff_scaling_tol}
\end{table}

\subsection{Simulation on multilevel grids}
\label{sec:multilevel}

\begin{figure}[tb!]
	\centering
	\includegraphics[width=0.45\textwidth]{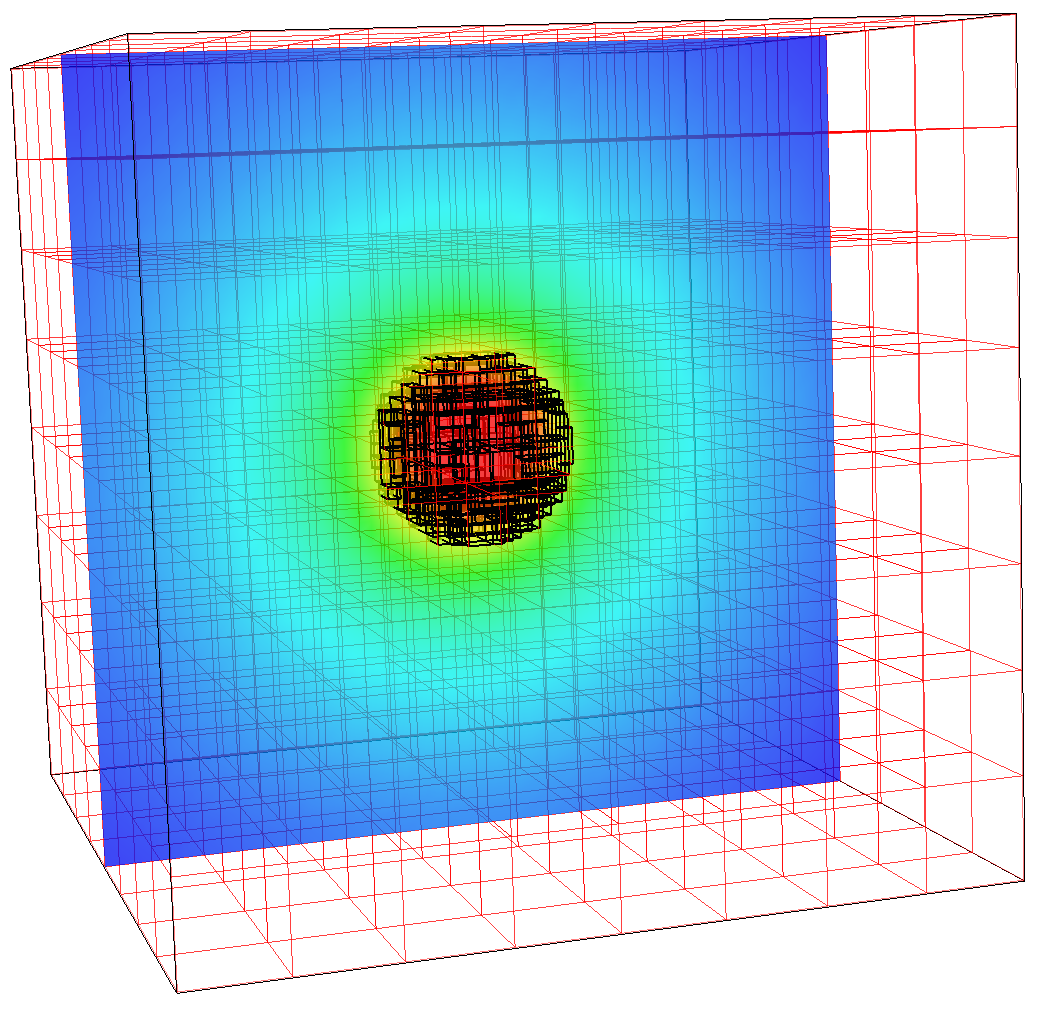}
	\caption{Distribution of boxes in 3D overlaid on a 2D slice of the methane mass fraction.  The cells within the black boxes have double the resolution of the cells within the red boxes.}
	\label{fig:multilevel}
\end{figure}

As a proof of concept and demonstration of capability, problem \eqref{eq:problem2} is also implemented in AMReX on multi-level grids.  The problem now uses a ball of pre-mixed methane fuel in the center of the domain, but is otherwise the same as from the single-level case.
Figure \ref{fig:multilevel} shows a 2D slice of the methane mass fraction at the initial time point, where we have overlaid the coarse and fine grid structures.
Unlike what is typically done with AMR, the multirate time stepping here is dictated by the accuracy requirements of the chemistry and diffusion and not by the CFL stability restriction.  The advection is computed explicitly, but the velocity is too small in this problem for CFL limitations to restrict step sizes.  The diffusion is computed implicitly and thus also does not limit the step size.  Therefore, the coarse step sizes for diffusion and advection and fine time step size for reactions are respectively the same for all levels of the multi-level grid.

To demonstrate the advantage of using both multirate time stepping and multi-level grids, versus without or with either option alone, we ran both ARK-IMEX and IMEX-MRI4-10 on single and two-level grids.  In the multilevel case, the fine grid only covers the fuel in the center, i.e., for a radius of $0.01$~cm extending from the center of the domain (which covers $[-0.05, 0.05]$ per side).  In the single level case, the domain contains $256^3$ grid cells divided into blocks of 32 cells on a side.  In the multilevel case, we use an effective $256^3$ resolution with a coarse domain containing $128^3$ grid cells, using blocks of 16 cells on a side to maintain the same number of MPI ranks as the single level case, along with another level of refinement.  A multilevel setup of our simulation is shown in Figure \ref{fig:multilevel}; we note that the fine grids cover roughly 5\% of the domain volume.  In the figure, the resolution is lower for clarity.  We leverage the composite multilevel MLABecLaplacian operator built into AMReX for implicit diffusion.  The solver tolerance is $10^{-12}$ for all cases.

Table \ref{tab:multilevel} shows that using either AMR or multirate alone takes considerably less time than single-rate time stepping on the single level grid.  However, combining the two techniques results in a considerably more efficient result than with either feature alone.

\begin{table}[tb!]
\begin{tabular}{l|c|r}
Integrator   & Total levels & Time \\ \hline
ARK-IMEX     & 1            & 8582 \\
ARK-IMEX     & 2            & 2363 \\
IMEX-MRI4-10 & 1            & 1507 \\
IMEX-MRI4-10 & 2            & 478
\end{tabular}
\caption{An illustration of the cost benefit of using both multirate in time and multilevel grids in space.
The effective resolution in each of these cases is 256$^3$ grid cells.
For the two-level cases, the fine grid covers only the fuel in the center as shown in Figure \ref{fig:multilevel}.
}
\label{tab:multilevel}
\end{table}

\section{Conclusions}
\label{sec:conclusions}

With the continued increases in compute power delivered by high performance systems, there continues to be increasing complexity in scientific simulations.  One area where this complexity manifests is in modeling systems with more physical processes that may operate at differing temporal and spatial scales.  Recent work in time integration methods has resulted in several new approaches to handling multiphysics systems with disparate time scales, and many of these are being developed into very efficient implementations within numerical libraries. In the meantime, work in mesh refinement technologies has resulted in highly efficient libraries for handling multiphysics systems with differing spatial scales. In this work, we evaluated new MRI and IMEX-MRI multirate time integration methods implemented within the SUNDIALS time integration library, along with block structured spatial mesh discretization approaches implemented within the BoxLib and AMReX libraries.  These evaluations were done within the context of high performance computing systems and also included comparison with spectral deferred correction methods and their multirate variants.
Results on two multiphysics test problems from combustion show that the new multirate methods match the accuracy of SDC and MRSDC methods while giving higher efficiency in most accuracy regimes.  In addition, results show that the multirate methods are more efficient than single rate methods for multiphysics systems posed within block structured mesh frameworks.  Lastly, the multirate and multi-spatial grid methods were shown to both improve overall solution efficiency, and their impacts were amplified when both approaches were simultaneously applied.

We expect these results to translate to other multiphysics systems with well separated time scales, such as are found in low Mach combustion systems and even in watershed models. In low Mach combustion, large time steps allow advection to be treated on a very slow time scale with implicit diffusion; reactions are even more challenging due to the wider separation of scales.
In watershed systems, hydrological models must be advanced with significantly larger time steps than land models.  Future work will consider application of multirate time integration methods to these and other applications.

Lastly, we note that the multirate methods can be made even more efficient through application of adaptive time step technologies, that can additionally relieve users from manual discovery of appropriate fast and slow step sizes for their applications.  These approaches have been very effective in single rate methods.  Current work is exploring how to bring such adaptivity into the multirate case.

\begin{acks}
This work was supported by the U.S. Department of Energy (DOE)
Office of Advanced Scientific Computing Research (ASCR)
via the Scientific Discovery through Advanced Computing (SciDAC) program FASTMath Institute.
Work at LLNL was performed under the auspices of the U.S. Department of Energy by Lawrence Livermore National Laboratory under contract DE-AC52-07NA27344, Lawrence Livermore National Security, LLC.
Work at LBL was performed under the auspices of the U.S. Department of Energy by Lawrence Berkeley National Laboratory under contract DE-AC02-05CH11231.

This document was prepared as an account of work sponsored by an agency of the United States government. Neither the United States government nor Lawrence Livermore National Security, LLC, nor any of their employees makes any warranty, expressed or implied, or assumes any legal liability or responsibility for the accuracy, completeness, or usefulness of any information, apparatus, product, or process disclosed, or represents that its use would not infringe privately owned rights. Reference herein to any specific commercial product, process, or service by trade name, trademark, manufacturer, or otherwise does not necessarily constitute or imply its endorsement, recommendation, or favoring by the United States government or Lawrence Livermore National Security, LLC. The views and opinions of authors expressed herein do not necessarily state or reflect those of the United States government or Lawrence Livermore National Security, LLC, and shall not be used for advertising or product endorsement purposes.
\end{acks}

\bibliographystyle{SageH}
\bibliography{AMReXSundials}

\end{document}